\documentclass[12pt]{article}
\usepackage[cp850]{inputenc}
\usepackage[T1]{fontenc}\usepackage{amsfonts}
\usepackage{graphicx}
\usepackage{amsmath}
\def\N{\mathbb{N}}\def\E{\mathbb{E}}\def\R{\mathbb{R}}

\author{G\'erard  Letac\thanks{Institut de Math\'ematiques de Toulouse, Universit\'e Paul Sabatier, 118 route de Narbonne 31062 Toulouse, \texttt{gerard.letac@math.univ-toulouse.fr}}}
\date{\today}
\title{Cumulative distribution functions for  the five simplest  natural exponential families}
\begin{document}  \maketitle

\begin{abstract} Suppose that the distribution of  $X_a$ belongs to a natural exponential family concentrated on the nonegative integers  and is such that $\E(z^{X_a})=f(az)/f(a)$. Assume  that $\Pr(X_a\leq k)$ has the form $c_k\int_a ^{\infty}u^k\mu(du)$ for some number $c_k$ and some positive measure $\mu,$ both independent of $a.$ We show that this asumption implies that the exponential family is either a binomial, or  the Poisson, or a negative binomial  family. Next, we study an analogous property for continuous distributions and we find that it is satisfied  if and only  the families are either Gaussian or Gamma. 
Ultimately, the proofs rely on the fact that only Moebius functions preserve the cross ratio, 


\end{abstract}

\section{Introduction}
\subsection{The three classical discrete examples}
If $X$ is valued in $\N=\{0,1,2, \ldots\}$ then there are three classical examples of an integral representation of the cumulative distribution function $\Pr(X\leq k).$ They are the following:

\begin{enumerate}
\item \textbf{Binomial case.} If $X\sim B(N,p)$ and if $a=p/(1-p)>0$ we
have
$$\E(z^X)=(1-p+pz)^N=\frac{(1+az)^N}{(1+a)^N}.$$ In particular for $k=0,\ldots,N$ we have 
\begin{eqnarray}\nonumber\Pr(X\leq k)&=&\frac{N!}{k!(N-k-1)!)}\int_p^1t^k(1-t)^{N-k-1}dt\\&=&\label{REP1}\frac{N!}{(N-k-1)!}\int_a^{\infty}\frac{1}{(1+u)^{N+1}}\frac{u^k}{k!}du.\end{eqnarray}

\item \textbf{Poisson case.} If $X\sim \mathcal{P}_{a}$ with $a>0$ we have 
$\E(z^X)=e^{a(z-1)}.$ In particular  for $k\in \N$ 
\begin{equation}\label{REP2}\Pr(X\leq k)=\int_{a}^{\infty}e^{-u}\frac{u^k}{k!}du\end{equation}

\item \textbf{Negative binomial case.}
If $X\sim NB(\lambda,a)$ with $0<a<1$ and $\lambda$ a positive real number, we have
$\E(z^X)=\frac{(1-a)^\lambda}{(1-az)^{\lambda}}.$ In particular for $k\in \N$ we have 
\begin{equation}\label{REP3}\Pr(X\leq k)=(\lambda)_{k+1}\int_a^1(1-u)^{\lambda-1}\frac{u^k}{k!}du\end{equation} where $(\lambda)_n=\lambda(\lambda+1)\cdots(\lambda+n-1)$ is the usual Pochhammer symbol. 

\end{enumerate}
These three statements \eqref{REP1},  \eqref{REP2}, \eqref{REP3} can be checked simply by taking the derivatives of  both sides with respect to $a.$ 
\subsection{Why the Poisson law and the binomial laws are so special?} One can be surprized by the simplicity of \eqref{REP1},  \eqref{REP2}, \eqref{REP3}.  In these three cases, we can observe the following: we start from a non zero power series $$f(z)=\sum_{k=0}^{\infty}p_kz^k$$ with $p_k\geq 0$ and a positive radius of convergence $R\in (0,\infty].$   We denote  $$A=\inf\{
k; p_k>0\} ,\ A+N=\sup\{
k; p_k>0\}\leq \infty.$$ We consider the natural exponential family $F_f=\{P_a; 0<a<R\}$ on $\N$ defined by 
$P_a(k)=p_ka^k/f(a)$ and we observe for each of these three examples that there exists a positive function $g(u)$ defined on $(0,R)$ and a sequence $(c_k ; A\leq k<A+N)$ 
such that, for $X\sim P_a$  for all $a\in (0,R)$ and for all $k $ such that $A\leq k<A+N ,$ we have

\begin{equation}\label{PROPERTY}\Pr(X\leq k)=c_k \int_a^Ru^kg(u)du.\end{equation}

\begin{itemize}
\item \textbf{Binomial case.} If $X\sim B(N,a/(1+a))$ then $A=0$, $R=\infty$ and 
$$c_k =N \binom{N-1}{k},\ g(u)=\frac{1}{(1+u)^{N+1}}.$$
\item \textbf{Poisson case.} If $X\sim \mathcal{P}_a$ then $A=0,$ $N=\infty$, $R=\infty$ and 
$$c_k =1/k!,\ g(u)=e^{-u}.$$

\item \textbf{Negative binomial case.} If $X\sim NB(\lambda,a)$ then $A=0,$ $N=\infty$,  $R=1$ and 
$$c_k = \frac{(\lambda)_{k+1}}{k!},\ g(u)=(1-u)^{\lambda-1}.$$

\end{itemize}

It is not correct to think that any distribution on $\N$ has this property \eqref{PROPERTY}. For instance if the law of $X$ belongs to the natural exponential family generated by the counting measure on $\{0,1,\ldots,N\}$ then  for $a\neq 1$ and $k\leq N$ we have
$$\Pr(X\leq k)=\frac{a^{k+1}-1}{a^{N+1}-1}$$ which does not satisfy \eqref{PROPERTY}.

One aim of this note is to show  that the three exponential families of distributions above are, up to translation,  the only ones which satisfy \eqref{PROPERTY}. A trivial remark about the integer $A$ is in order: suppose that 
\eqref{PROPERTY} holds with $A>0.$ Then trivially $X-A$ will satisfy \eqref{PROPERTY} with $g_A(u)=u^Ag(u).$ Therefore  we can assume $A=0$, that is $p_0>0,$ in the study of the property \eqref{PROPERTY}. More specifically, replacing the positive density $g$ in 
\eqref{PROPERTY} by  an arbitrary positive measure $\mu$ one obtains the following result.

\vspace{4mm}\noindent \textbf{Theorem 1.} Let $f,$ $N\leq \infty$,  $R$ and the corresponding natural exponential family   $F_f=\{P_a; 0<a<R\} $ be as before.  Assume that $p_0=f(0)>0$ and that $f$ is not constant. Suppose that there exists a positive measure $\mu(du)$ on $(0,R)$ and a sequence $(c_k)_{0\leq k<N}$ such that for all $k<N$, for all $a\in (0,R)$ and $X\sim P_a$ one has 
\begin{equation}\label{PROPERTY2}\Pr(X\leq k)=c_k\int_{(a,R)}u^k\mu(du)\end{equation}
 Then the  exponential family $F_f$ is either binomial, or Poisson, or negative binomial. 
 
\subsection{The continuous analogues} We describe now the second result of this note.  Let us  consider a property similar to \eqref{PROPERTY} or \eqref{PROPERTY2}. Let $\nu$ be a positive measure on $\R$ not concentrated on one point, with Laplace transform 
$$L_{\nu}(\theta)=\int_{\R}e^{\theta x}\nu(dx)\leq \infty,$$ and such that the open interval $(\alpha,\beta)$ which is the interior of the set 
$\{\theta; L_{\nu}(\theta)<\infty\}$ is not empty. The measure $\nu$ is not necessarily bounded.  Denote $\kappa_ {\nu}(\theta)=\log L_{\nu}(\theta).$ For $\alpha<\theta<\beta$ consider the probability
$$P_{\theta}(dx)=e^{\theta x-\kappa_ {\nu}(\theta)}\nu(dx)$$ and the natural  exponential family $$F(\nu)=\{P_{\theta}\ ;\  \alpha<\theta<\beta\}.$$ Suppose now that there exists a function $c$ on $\R$ and a positive measure $\mu(du)$ on  $(\alpha,\beta)$ such that for all $x\in \R$, for all $\theta$ in $(\alpha,\beta)$ and for $X\sim P_{\theta}$ we have 
\begin{equation}\label{PROPERTY3}
\Pr(X\leq x)=c(x)\int_{(\theta,\beta)}e^{ux}\mu(du).
\end{equation} Note that here we do not assume that $\nu$ has no atoms. This theorem is called a continuous analog of Theorem 1 only because  \eqref{PROPERTY3} has to be true for all real $x.$ 
Here are the two classical examples
\begin{enumerate}
\item\textbf{The Gaussian case.} We take $\nu(dx)=e^{-\frac{x^2}{2\sigma^2}}\frac{dx}{\sigma\sqrt{2\pi}}.$ Then $(\alpha,\beta)=\R$ and 
$$P_{\theta}(dx)=e^{-\frac{1}{2\sigma^2}(x-\theta \sigma^2)^2}\frac{dx}{\sigma\sqrt{2\pi}}.$$ As a consequence, for $X\sim P_{\theta}$ we have
$$\Pr(X<x)=e^{-\frac{x^2}{2\sigma^2}}\int_{\theta}^{\infty}e^{ux}e^{-\frac{u^2\sigma^2}{2}}\frac{du}{\sigma\sqrt{2\pi}}.$$ 
Therefore $c(x)=e^{-\frac{x^2}{2\sigma^2}}$ and $\mu(du)=e^{-\frac{u^2\sigma^2}{2}}\frac{du}{\sigma\sqrt{2\pi}}.$

\item\textbf{The Gamma case.} For some $p>0$ we take $\nu(dx)=x^{p-1}1_{(0,\infty)}(x)dx/\Gamma(p).$Then $(\alpha,\beta)=(-\infty,0)$ and 
$$P_{\theta}(dx)=e^{\theta x}(-\theta)^p x^{p-1}\frac{dx}{\Gamma(p)}.$$ As a consequence, for $X\sim P_{\theta}$ we have for $x>0$
$$\Pr(X<x)=\frac{x^{p-1}}{\Gamma(p)}\int_{\theta}^0e^{ux}(-u)^{p-1}du$$ 
Therefore $c(x)=x^{p-1}1_{(0,\infty)}(x)$ and $\mu(du)=(-u)^{p-1}du.$

\end{enumerate}

Again not all natural exponential families have property \eqref{PROPERTY3}; for instance if $\nu$ is the uniform distribution on $( 0,1)$, we see that it is impossible to find $c$ and $\mu$ such that for all $\theta$ in $\R$ and all $x\in (0,1)$ such that
$$ \frac{e^{\theta x}-1}{e^{\theta}-1}=c(x)\int_{\theta}^{\infty}e^{ux}\mu(du).$$

These two examples where \eqref{PROPERTY3} holds happen to be the only possible ones, up to translation. More specifically:

\vspace{4mm}\noindent \textbf{Theorem 2.} Let $\nu$, $(\alpha,\beta)$  and the corresponding natural exponential family   $F(\nu)$ be as before.   Suppose that there exists a positive measure $\mu(du)$ on $(\alpha,\beta)$ and a function $c$ on $\R$ such that for all $x\in R$,  for all $\theta\in (\alpha,\beta)$ and $X\sim P_{\theta}$  equality \eqref{PROPERTY3} holds.  
 Then up to translation, the  exponential family $F(\nu)$ is either Gaussian or Gamma.

\vspace{4mm}\noindent To prove these theorems in Sections 2 and 3, we  need the following  lemma, that will be shown in Section 4. Some comments are given in Section 5.

\vspace{4mm}\noindent \textbf{Lemma.} Let $K,A,B$ be three   functions defined on an open interval $I$ such that for all $v< u$ with $u,v\in I$ we have 
\begin{equation}\label{CROSS}\frac{K(u)-K(v)}{u-v}=A(u)B(v).\end{equation} Assume also that $A$ and $B$ are continuous and strictly positive on  $I.$ Then $K$ is a Moebius function, that is  of the form $K(u)=\frac{au+b}{cu+d}$ where $ad-cd\neq 0.$
\section{Proof of Theorem 1.} 
For fixed $z\in (0,1)$ we observe that for all $u\in (0,R).$
\begin{equation}\label{DERIVEE}\frac{d}{du}\frac{f(uz)}{f(u)}=\frac{zf'(uz)f(u)-f(uz)f'(u)}{f(u)^2}.\end{equation}
 The function $u\mapsto f(uz)/f(u)$ is strictly decreasing. To see this, recall from the properties of the natural exponential families that the cumulant function $\kappa$ defined on $(-\infty,\log R)$ by $\kappa(\log u)=\log f(u)$ is 
strictly convex since $f$ is not constant and since $f(0)>0$. Therefore $$ \kappa'(\log u+\log z)-\kappa'(
\log u)<0,$$ which  implies that $\frac{d}{du}\frac{f(uz)}{f(u)}<0$. As a consequence $\lim_{u\uparrow R}\frac{f(uz)}{f(u)}$ always exists and is less than 1. 
Therefore for $a\in (0,R)$ we have
\begin{equation}\label{DERIVEE1}
\frac{f(az)}{f(a)}=\lim_{u\uparrow R}\frac{f(uz)}{f(u)}-\int_a^R\frac{d}{du}\frac{f(uz)}{f(u)}du.
\end{equation}
For clarity now we distinguish the cases $N$ finite and $N$ infinite.

\subsection{The $N$ finite case.} Obviously $R=\infty
.$ Note that here $\lim_{u\uparrow R}\frac{f(uz)}{f(u)}=z^N.$ Applying \eqref{DERIVEE1} and using $$ \sum_{k=0}^{\infty}\Pr(X\leq k)z^k=\frac{1}{1-z}\E(z^X)$$ we get

$$\sum_{k=0}^{\infty}\Pr(X\leq k)z^k=\frac{z^N}{1-z}-\frac{1}{1-z}\int_a^{\infty}\frac{d}{du}\frac{f(uz)}{f(u)}du$$ which implies

\begin{equation}\label{FINI}\sum_{k=0}^{N-1}\Pr(X\leq k)z^k=-\frac{1}{1-z}\int_a^{\infty}\frac{d}{du}\frac{f(uz)}{f(u)}du\end{equation} Now let us introduce the polynomial 
$C(u)=\sum_{k=0}^{N-1}c_ku^k$ and let us apply the hypothesis \eqref{PROPERTY2} to \eqref{FINI}. For all $a>0$ we obtain
$$\int_{(a,\infty)} C(uz)\mu(du)=-\frac{1}{1-z}\int_a^{\infty}\frac{d}{du}\frac{f(uz)}{f(u)}du.$$ Taking derivative with respect to $a$ and using \eqref{DERIVEE} we get our fundamental equation 

\begin{equation}\label{FUND}C(uz)\mu(du)=-\frac{zf'(uz)f(u)-f(uz)f'(u)}{(1-z)f(u)^2}du.\end{equation}
Equality \eqref{FUND} shows that actually $\mu(du)=g(u)du$ for some analytic  function $g$. Since we have seen that $\frac{d}{du}\frac{f(uz)}{f(u)}<0,$ this finally proves that for all $u\in (0,R)$ we have the important fact $g(u)>0.$

\subsection{The $N$ infinite case.} We are now prove  equation \eqref{FUND} when $N=\infty$  in a quite similar way as $N<\infty.$  We introduce the sum 
$$C_n(u)=\sum_{k=0}^nc_ku^k.$$ Since $c_k\geq 0$ then for $u\in (0,R)$ the function $C(u)=\lim_nC_n(u)$ exists in $(0,\infty].$
Since from \eqref{PROPERTY2} we can write
$$\sum_{k=0}^n\Pr(X\leq k)u^k=\int_a^RC_n(uz)\mu(du)$$ we can pass to the limit when $n\to \infty$ and we obtain by monotone convergence 
\begin{equation}\label{INFINI}\sum_{k=0}^{\infty}\Pr(X\leq k)u^k=\int_a^RC(uz)\mu(du)\end{equation}
where both sides are finite since $0<z<1.$ We use now

\begin{equation}\label{INFINI1}\sum_{k=0}^{\infty}\Pr(X\leq k)u^k=\frac{1}{1-z}\lim_{u\uparrow R}\frac{f(uz)}{f(u)}-\frac{1}{1-z}\int_a^R\frac{d}{du}\frac{f(uz)}{f(u)}du.\end{equation}
Next, let us compare \eqref{INFINI} and  \eqref{INFINI1} and let us take derivatives with respect to $a.$ Equation  \eqref{FUND} for $N=\infty$ follows. A similar reasoning as in the case $N$ finite shows that $\mu(du)=g(u)du$ for some strictly positive and analytic function $g.$  Furthermore, from \eqref{FUND} we get that $C(u)$ is finite on $(0,R).$

\subsection{Solution of the functional equation \eqref{FUND} .}  In  \eqref{FUND}, where of course we replace $\mu(du)$ by $g(u)du,$ we denote $v=uz$: thus $0<v<u<R.$ Since $f(u)$ is positive we divide both sides of  \eqref{FUND}  by $f(u)f(v)$ and we  denote
$$A(u)=g(u)f(u),\ B(v)=\frac{C(v)}{f(v)},\ K(u)=\frac{uf'(u)}{f(u)}.$$ Note that $A(u)>0,\ B(v)>0$ for all $u,v$ in $(0,R)$ and that $A$ and $B$ are continuous.

With these notations \eqref{FUND} becomes 
$$A(u)B(v)=\frac{K(u)-K(v)}{u-v}.$$ From the lemma applied to $I=(0,R)$ we know that  there exists $a,b,c,d$ such that $K(u)=(au+b)/cu+d).$ and $ad-bc\neq 0.$ We have
$$\frac{f'(u)}{f(u)}=\frac{au+b}{u(cu+d)}.$$ We discuss the various  particular cases:
\begin{itemize}
\item \textbf{ $d=0.$} This leads to $f(u)=Au^Be^{C/u}$ with $C\neq 0.$ Such an $f$ is not analytic around $0$ and this case is excluded.

\item \textbf{ $c=0.$} This leads to $f(u)=Au^Be^{Cu}$ with $C\neq 0.$ Since $f(0)>0$ we have $B=0$ and we are in the Poisson case. 
 \item \textbf{ $c\neq 0, \ d\neq 0.$} This leads to $f(u)=Au^B(C\pm u)^D.$ With  $A=f(0)>0$ we have $B=0.$ The  fact that $f$ must have a power series expansion with non negartive coefficients implies that either $f(u)=A(C+u)^D$ where $C>0$ and $D$ is a positive integer (and we are in the binomial case) or 
$f(u)=A(C+u)^D$ where $C>0$ and $D$ is a negative real number  (and we are in the negative  binomial case) .$\square$
\end{itemize}
\section{Proof of Theorem 2} We multiply both sides of \eqref{PROPERTY3} by $e^{zx}$ where $\alpha-\theta<z<0$  and we integrate in $x$ on the whole $\R$. The left hand side of  \eqref{PROPERTY3} is transformed as follows:
\begin{eqnarray} 
\nonumber\int_{\R}e^{xz}\Pr(X\leq x)dx&=&\int_{\R} e^{xz}\left(\int_{(-\infty,x]}e^{\theta t-\kappa_ {\nu}(\theta)}\nu(dt)\right)dx\\&=&\label{LEFT}e^{-\kappa_ {\nu}(\theta)}\int_{\R}e^{\theta t}\left(\int_{[t,\infty)}e^{xz}dx\right)\nu(dt)=\frac{1}{-z}e^{\kappa_ {\nu}(\theta+z)-\kappa_ {\nu}(\theta)}\end{eqnarray}The  right hand side of  \eqref{PROPERTY3} is similarly transformed:
\begin{eqnarray} 
\nonumber\int_{\R}c(x)e^{zx}\left(\int_{(\theta,\beta)}e^{ux}\mu(du)\right)dx&=&\int_{(\theta,\beta)}\left(\int_{\R}e^{u+z}c(x)dx\right)\mu(du)\\&=&\label{RIGHT}\int_{(\theta,\beta)}e^{\kappa_ c(u+z)}\mu(du)
\end{eqnarray}
Comparing \eqref{LEFT} and \eqref{RIGHT} and differentiating in $\theta$ we get
\begin{equation}\label{FUND2}
\frac{1}{-z}e^{\kappa_ {\nu}(\theta+z)-\kappa_ {\nu}(\theta)}(\kappa' _{\nu}(\theta+z)-\kappa' _{\nu}(\theta))d\theta=-e^{\kappa_ c(\theta+z)}\mu(d \theta)\end{equation} which shows that $\mu$ is absolutely continuous. We therefore denote $g(\theta)d\theta=\mu(d \theta)$  and \eqref{FUND2} becomes 
\begin{equation}\label{FUND3}
\frac{1}{-z}e^{\kappa_ {\nu}(\theta+z)-\kappa_ {\nu}(\theta)}(\kappa' _{\nu}(\theta+z)-\kappa' _{\nu}(\theta))=-e^{\kappa_ c(\theta+z)}g(\theta)\end{equation} 
Since $\nu$ is not concentrated on one point, the function $\kappa_ {\nu}$ is strictly convex on 
$(\alpha,\beta)$ and therefore $\kappa' _{\nu}(\theta+z)-\kappa' _{\nu}(\theta)<0.$ As a consequence \eqref{FUND3} shows that $g$ is continuous and strictly positive on $(\alpha,\beta).$ In order to use the lemma, we modify \eqref{FUND3} by using the  notation $\theta=u,\ \theta+z=v,$ implying $$ \alpha<v<u<\beta.$$ We obtain  
\begin{equation}\label{FUND4}
\frac{1}{u-v}e^{\kappa_ {\nu}(v)-\kappa_ {\nu}(u)}(\kappa' _{\nu}(u)-\kappa' _{\nu}(v)=e^{\kappa_ c(v)}g(u)\end{equation} 
We indeed apply the lemma to \eqref{FUND4} with the notations 
$$I=(\alpha,\beta),\ A(u)=g(u)e^{\kappa_ {\nu}(u)},\ B(v)=e^{\kappa_ c(v)-\kappa_ c(v)},\ K(u)=\kappa' _{\nu}(u)$$
and we can claim that $\kappa' _{\nu}(u)$ is a Moebius function restricted to $(\alpha,\beta).$  A standard discussion leads to the result: either $\kappa' _{\nu}$ is affine and we are in the Gaussian case, or not. In this second case we have 
$\kappa' _{\nu}(\theta)=A-\frac{p}{\theta-B}$ where $p>0$ to insure convexity. Consider first  the case $B \leq \alpha.$ It leads to  $$L_{\nu}(\theta)=e^{A\theta}\frac{C}{(\theta-B)^p}$$ which is the case of a Gamma family with shape parameter $p$ and translated by $A.$ The case $B\geq \beta$ leads to  $L_{\nu}(\theta)=e^{A\theta}\frac{C}{(B-\theta)^p}$ corresponding to the opposite of a Gamma distribution. Since its support is $(-\infty,A],$ such a family  cannot satisfy \eqref{PROPERTY3} for $x>A$ and this case has to be excluded. $\square$

\section{Proof of the lemma.}  Before giving  the proof, let us remark  that the  lemma would be much easier  to show if \eqref{CROSS} was supposed to be true for all $u,v$ in $I$. Indeed  \eqref{CROSS} without the restriction  $v<u $ would imply that $K$ preserves the cross ratio of a sequence of four real distinct numbers $(x_1,x_2,x_3,x_4)$ defined as 
$$[x_1,x_2,x_3,x_4)]=\frac{x_3-x_1}{x_3-x_2}\times \frac{x_4-x_2}{x_4-x_1}$$ since
\begin{equation}\label{MOEBIUS}\frac{[K(u),K(x_2),K(x_3),K(x_4)]}{[u,x_2,x_3,x_4)]}=\frac{A(x_3)B(u)A(x_4)B(x_2)}{A(x_4)B(u)A(x_3)B(x_2)}=1.\end{equation} The lemma would be then proved by using the simple fact that a function $K$ which satisfies \eqref{MOEBIUS} is necessarily a Moebius function.

\vspace{4mm}\noindent
We now prove the lemma in direr circumstances  with the  less restrictive hypothesis that \eqref{CROSS} is true at least for $v<u$.  Because of the continuity asumptions on $A$ and $B$ clearly $K'$ exists and $ K'(u)=A(u)B(u).$  Applying $\frac{\partial ^2}{\partial u\partial v}$ to both sides of the equality \begin{equation}\label{CROSS1}K(u)-K(v)=(u-v)A(u)B(v)\end{equation} on the set 
$D=\{(u,v)\in I^2: v\leq u\}$ we obtain 
\begin{equation}\label{CROSS2}A'(u)B(v)-A(u)B'(v)=(u-v)A'(u)B'(v)\end{equation}

Suppose that there exists $v_0$ such that $B'(v_0)=0.$ Then from \eqref{CROSS2} we have that $A'(u)=0$ and $A(u)=A(v_0)$ for all $u\geq v_0$ and therefore, from \eqref{CROSS1}, $K$ is affine on $u\geq v_0,$ with $K'(u)=A(v_0)B(v)$ for all $u>v_0$ and all $v<u.$ This implies that $B$ is a constant on $I$, and $A$ is a constant as well. Finally $K$ is affine on $I$, a particular case of Moebius function.

Suppose now that there exists $u_0$ such that $A'(u_0)=0.$ The same reasoning shows that $K$ is still an affine function.

From now on we may assume that $A'(u)B'(v)\neq 0$ on $D.$ Dividing both sides of \eqref{CROSS2}
 by $A'(u)B'(v)$ one obtains on $D$ that

\begin{equation}\label{CROSS3}\frac{A(u)}{A'(u)}+u=\frac{B(v)}{B'(v)}+v\end{equation}
We are in position to apply the principle of separation of variables and we can claim that there exists a constant $d$ such that  both sides of \eqref{CROSS3} are equal to $-d$ for all $u$ and $v$ in $I$ without the restriction $(u,v)\in D.$ From this we get easily the existence of constants $A$ and $B$ such that
$$A(u)=\frac{A}{u+d},\ B(v)=\frac{B}{v+d},\  K'(u)=A(u)B(u)=\frac{AB}{(u+d)^2}$$ implying that $K$ is a Moebius function. $\square$

\section{Comments} Many  introductory textbooks in probability and statistics, as well as Wikipedia,  are mentioning \eqref{REP1}  and \eqref{REP2}, while  \eqref{REP3}  appears less frequently.  A reference is Abramowitz and Stegun (1965), formulas 26.5.24 and 26.5.26 page 945. The similarity of these formulas \eqref{REP1},\eqref{REP2},\eqref{REP3}   is not really commented. Sometimes it is said that \eqref{REP1},  \eqref{REP2} and  \eqref{REP3} use respectively  the 'truncated beta distribution' of the second kind with  parameters $(k+1,N-k)$, the 'truncated gamma distribution' of parameter $k+1$ and  the 'truncated beta distribution' of the first kind with  parameters $(k+1,\lambda).$ These three continuous distributions are qualified of 'conjugate' with respect to the binomial, Poisson and negative binomial distributions respectively, while the definition  of conjugacy is not given. Sometimes it is observed that if $\lambda$ and $k$ are positive integers, if $0<a<1$ and if $X\sim B(k+\lambda,1-a)$ , $Y\sim NB(\lambda,a)$ then 
$$\Pr(X\geq \lambda)=\Pr(Y\leq k).$$ This formula  can be checked with \eqref{REP1} and  \eqref{REP3} and lead some authors to say in an unspecified sense  that the binomial and negative binomial families are 'inverse' of each other. Also, $\eqref{REP1}$ is related to the $k$ th statistic $U^{(N)}_k$ of a sample $(U_1,\ldots,U_N)$ of iid rv which are uniform on $(0,1)$, by mean of the formula 
$$\Pr(U^{(N)}_k>p)=\Pr(X\leq k)$$
when $X\sim B(N,p).$

Finally, there are  other common characterizations of these three discrete exponential families. One is rather trivial with $\Pr(X=k+1)/\Pr(X=k)=(Ak+B) /(k+1).$ Another one comes  from the point of view of the variance functions, in the sense of Carl Morris (1982). Indeed, these variance functions are respectively 
$V(m)=m-\frac{m^2}{N}$ on the mean domain $(0,N)$, $V (m)=m$ and $V(m)=m+\frac{m^2}{\lambda}$ on the mean domain $(0,\infty).$ Therefore, they have the common form $V(m)=m+cm^2.$ We have not been able to link the characterization of the present note with the characterization by variance functions.

We have completed the characterization of Theorem 1 with its continuous analog in Theorem 2, with the less striking property \eqref{PROPERTY3}.  It would have been possible but cumbersome to gather the two theorems in one. It was a surprise to have to use the lemma again in the continuous case.

Normal and Gamma exponential families also belong to the Morris families, the ones which have quadratic variance functions,  respectively $V(m)=\sigma^2$ for the Gauss case and $V(m)=\frac{1}{p}(m-A)^2$ on the mean domain  $(A,\infty)$ for the Gamma case. A confirmation of the guess that the variance functions are not the good tool to obtain the characterizations ot Theorems 1 and 2 is the fact that the 6th family of Morris, the one with variance function of the form $p+\frac{m^2}{p}$ does not appear here.

\section {References}
\begin{itemize}
\item \textsc{Abramowitz, M. and Stegun, I.} (1965) \textit{Handbook of Mathemarical Functions}, Dover, New York.
\item  \textsc{Morris, C. N.} (1982) Natural exponential families with quadratic variance functions. \textit{Ann. Statist.}, \textbf{10} 65-80.

\end{itemize}

\end{document}